\documentclass[12pt, journal, draftclsnofoot, onecolumn]{IEEEtran}
\ifCLASSOPTIONcompsoc
\else
\fi

\ifCLASSINFOpdf
\else
\fi
\usepackage[parfill]{parskip}    
\usepackage{amsmath,amssymb,latexsym,enumerate, mathrsfs}
\usepackage{graphicx}
\usepackage{graphicx}
\usepackage{array}
\usepackage{times}
\usepackage{hyperref}
\newtheorem{theorem}{\bf Theorem}
\newtheorem{lemma}{\bf Lemma}
\newtheorem{corollary}{\bf Corollary}

\newtheorem{remark}{Remark}
\newtheorem{definition}{Definition}

\def\etal{\mbox{et al.}}

\begin{document}
%
\title{On a connection between the reliability of multi-channel systems and the notion of controlled-invariance entropy}

\author{Getachew~K.~Befekadu~\IEEEmembership{}
\IEEEcompsocitemizethanks{\IEEEcompsocthanksitem G. K. Befekadu is with the Department
of Electrical Engineering, University of Notre Dame, Notre Dame, IN 46556, USA.\protect\\
E-mail: gbefekadu1@nd.edu\protect\\
Version - February 25, 2013.}
\thanks{}}

\markboth{}%
{Shell \MakeLowercase{\textit{et al.}}: Bare Advanced Demo of IEEEtran.cls for Journals}
\IEEEcompsoctitleabstractindextext{%
\begin{abstract}
The purpose of this note is to establish a connection between the problem of reliability (when there is an intermittent control-input channel failure that may occur between actuators, controllers and/or sensors in the system) and the notion of controlled-invariance entropy of a multi-channel system (with respect to a subset of control-input channels and/or a class of control functions). We remark that such a connection could be used for assessing the reliability (or the vulnerability) of the system, when some of these control-input channels are compromised with an external ``malicious" agent that may try to prevent the system from achieving more of its goal (such as from attaining invariance of a given compact state and/or output subspace).
\end{abstract}

}

\maketitle

\IEEEdisplaynotcompsoctitleabstractindextext

%
\IEEEpeerreviewmaketitle

\section{Introduction} \label{sec:Intro}
In the last decades, the notions of measure-theoretic entropy and topological entropy have been intensively studied in the context of measure-preserving transformations or continuous maps (e.g., see \cite{Wal82}, \cite{Sin94} and \cite{Dow11} for the review of entropy in ergodic theory as well as in dynamical systems). For instance, Adler \etal\, (in the paper \cite{AdlKM65}) introduced the notion of topological entropy as a topologically invariant conjugacy, which is an analogue to the notion of measure-theoretic entropy, for measuring the rate at which a continuous map in a compact topological space generates initial-state information. Later, \cite{Din70} and \cite{Bow71} gave a new, but equivalent, definition of topological entropy for continuous maps that led to proofs for connecting the topological entropy with that of measure-theoretic entropy. 

In the recent paper \cite{ColKa11}, the authors have introduced the notion of invariance entropy for continuous-time systems as a measure of information that is necessary to achieve invariance of a given state (or output) subspace (i.e., a measure of how open-loop control functions have to be updated in order to achieve invariance of a given subspace of the state space). In the present paper, we explore this concept, which is closely related to the notion of topological feedback entropy (see \cite{NaiEMM04} and \cite{HagN12}), for assessing the reliability of a multi-channel system when there is an intermittent control-channel failure that may occur between actuators, controllers and/or sensors in the system. Specifically, we provide conditions on the minimum rate at which the multi-channel system can generate information with respect to a subset of control-input channels and/or a class of control functions when the system states are restricted to a given controlled-invariant subspaces. Here, it is important to note that the intermittent control-channel failures may not necessarily represent any physical failures within the system. Rather, this can also be interpreted as external ``malicious" agent who is trying to prevent the system from achieving more of its goal, i.e., from attaining invariance of a given state (or output) subspace.

With the emergence of networked control systems (e.g., see \cite{AntBa07}), these notions of entropy have found renewed interest in the research community (e.g., see \cite{NaiEMM04}, \cite{Sav06} and \cite{ColKa09}). Notably, in the paper \cite{NaiEMM04}, Nair \etal\, have introduced the notion of topological feedback entropy, which is based on the ideas of \cite{AdlKM65}, to quantify the minimum rate at which deterministic discrete-time dynamical systems generate information relevant to the control objective of set-invariance. More recently, the notion of controlled-invariance entropy (as well as the notion of almost invariance entropy) has been studied for continuous-time control systems in \cite{ColKa09}, \cite{ColHe11}, \cite{Kaw11a} and \cite{ColKa11} based on the metric-space technique of \cite{Bow71}. It is noted that such an invariant entropy provides a measure of the smallest growth rate for the number of open-loop control functions that are needed to confine the states within an arbitrarily small distance (in the sense of gap metric) from a given subspace. For discontinuous systems, we also note that the notion of topological entropy has been investigated with respect to piecewise continuous piecewise monotone transformations (e.g., see \cite{Kop05} and \cite{MisZi92}).

The paper is organized as follows. In Section~\ref{sec:Prelim}, we present preliminary results on the invariance entropy of multi-channel systems with respect to a set of control-input channels and a class of control functions. Section~\ref{sec:main} presents the main results -- where we provide conditions on the information this is necessary for achieving invariance of the multi-channel system states in (or near) a given subspace.

\section{Preliminaries} \label{sec:Prelim}
\subsection{Notation} 
For $A \in \mathbb{R}^{n \times n}$, $B \in \mathbb{R}^{n \times r}$ and a linear space $\mathscr{X}$, the supremal $(A, B)$-invariant subspace contained in $\mathscr{X}$ is denoted by $\mathscr{V}^*\triangleq\sup \mathscr{V} \left(A, B; \mathscr{X} \right)$. For a subspace $\mathscr{V} \subset\mathscr{X}$, we use $\langle A \,|\, \mathscr{V} \rangle$ to denote the smallest invariant subspace containing $\mathscr{V}$.

\subsection{Problem formulation}
Consider the following generalized multi-channel system  
\begin{align}
 \dot{x}(t) &= A x(t) + \sum_{j \in \mathcal{N}} B_j u_j(t), \quad x(t_0)=x_0, \quad t \in [t_0, +\infty ), \label{Eq1}
\end{align}
where $A \in \mathbb{R}^{n \times n}$, $B_j  \in \mathbb{R}^{n \times r_j}$, $x(t) \in \mathscr{X} \subset \mathbb{R}^{n}$ is the state of the system, $u_j(t) \in \mathscr{U}_j \subset \mathbb{R}^{r_j}$ is the control input to the $j$th-channel and $\mathcal{N} \triangleq \{1, 2, \ldots, N\}$ represents the set of controllers (or the set of control-input channels) in the system.

Let us introduce the following class of admissible controls that will be used in the sequel
\begin{align}
 \mathscr{U} \subseteq \biggm \{ u \in \prod_{i \in \mathcal{N}} L_{\infty}(\mathbb{R},\, \mathbb{R}^{r_i}) \, \biggm \vert \, u_{\neg j}(t) \in \mathscr{U}_{\neg j} \triangleq \prod_{i \in \mathcal{N}_{\neg j}} \mathscr{U}_i ~~ \text{{\em for almost all}}& \,\, t \in [0,\, \infty) ~\, \text{and}\notag\\
                                                                                                                   &\forall j \in \mathcal{N} \cup \{0\} \biggm \}, \label{Eq2}
\end{align}
where $u_{\neg 0}(t)=\bigl(u_1(t), \, u_{2}(t), \, \ldots \,u_N(t)\bigr)$ and $u_{\neg i}(t)=\bigl(u_1(t), \, \dots \, u_{i-1}(t), \, u_{i+1}(t), \, \ldots \,u_N(t)\bigr)$ for $i \in \mathcal{N}$. Moreover, $\mathcal{N}_{\neg 0} \triangleq \mathcal{N}$ and $\mathcal{N}_{\neg j} \triangleq \mathcal{N}\setminus\{j\}$ for $j=1, 2, \ldots, N$.

In the remainder of this subsection, we provide some results from geometric control theory (e.g., see \cite{BasM69}, \cite{Won79}, \cite{BasM92} and \cite{TreStH01} for details about this theory).

\begin{definition}
Let $\mathscr{V}_j \subset \mathscr{X}$ for $j \in \mathcal{N} \cup\{0\}$.
\begin{enumerate} [(i)]
\item If $\mathscr{V}_j$ is $(A)$\,-\,invariant, then $A \mathscr{V}_j \subset \mathscr{V}_j$.
\item If $\mathscr{V}_j$ is $(A, \, B_{\neg j})$-invariant, then $A \mathscr{V}_j \subset \mathscr{V}_j +\mathscr{B}_{\neg j}$, where $B_{\neg 0}\triangleq\bigl[\begin{array}{cccc} B_1 & B_{2} & \ldots & B_N \end{array} \bigr]$, $B_{\neg j}\triangleq~\bigl[\begin{array}{ccc} B_1\,\ldots\,B_{j-1} & B_{j+1} \,\ldots\,B_N \end{array} \bigr]$ and \, $\mathscr{B}_{\neg j} \triangleq \operatorname{Im}B_{\neg j}$ for $j \in \mathcal{N}$.\footnote{In this paper, we consider a case in which one of the controllers is extracted due to an intermittent failure. However, following the same discussion, we can also consider when the fault is associated with at most two or more possible controllers in the system.}
\end{enumerate}
\end{definition}

The following lemma, which is a well-known result, will be stated without proof (e.g., see \cite{Won79} or \cite{BasM92}).
\begin{lemma} \label{L1}
Suppose $\mathfrak{I}_j \left(A, B_{\neg j}; \mathscr{F} \right)$ is a family of $(A, B_{\neg j})$-invariant subspaces for $j \in \mathcal{N}$. Then, every subspace $\mathscr{F} \subset \mathscr{X}$ contains a unique supremal $(A, B_{\neg j})$-invariant subspace which is given by $\mathscr{V}_j^*=\sup \mathfrak{I}_j \left(A, B_{\neg j}; \mathscr{F} \right)$ for each $j \in \mathcal{N}$.
\end{lemma}

Then, we state the following result which is a direct application of Lemma~\ref{L1}.
\begin{theorem} \label{T1}
Let $\mathscr{V}_j \subset \mathscr{X}$ for each $j \in \mathcal{N}$. Then, $\mathscr{V}_j$ is a member of the subspace families $\mathfrak{I}_j(A, B_{\neg j}; \mathscr{X})$ that preserves the property of $(A, B_{\neg j})$-invariant (i.e., $\mathscr{V}_j \in \mathfrak{I}_j (A, B_{\neg j}; \mathscr{X})$), if and only if
\begin{align}
A \mathscr{V}_j \subset \mathscr{V}_j + \mathscr{B}_{\neg j}, ~~ \forall j \in \mathcal{N}. \label{Eq3} 
\end{align}
\end{theorem}
\begin{IEEEproof}
The proof follows the same lines of argument as that of  Wonham (see \cite{Won79} p. 88). Suppose $\mathscr{V}_j \in \mathfrak{I}_j (A, B_{\neg j}; \mathscr{X})$ and let $ v^j \in \mathscr{V}_j$ for $j \in \mathcal{N}$, then $(A, \sum_{i \in \mathcal{N}_{\neg j}} B_i K_i)v^j = w^j$ for some $w^j \in \mathscr{V}_j$, i.e.,
\begin{align}
Av^j = w^j -  \sum_{i \in \mathcal{N}_{\neg j}} B_i K_i v^j \in \mathscr{V}_j + \mathscr{B}_{\neg j}, \label{Eq4}
\end{align}
On the other hand, let $\{v_1^j, v_2^j, \ldots, v_\mu^j \}$ be a basis for $\mathscr{V}_j$ for $j=1, 2, \ldots, N$. Suppose that \eqref{Eq3} holds true. Then, there exist $w_k^j \in \mathscr{V}_j$ and $u_k^{\neg j} \in \mathscr{U}_{\neg j}$ for $k \in \{1, 2, \dots, \mu^j \}$ such that
\begin{align}
Av_k^j = w_k^j -  B_{\neg j} u_k^{\neg j}, \qquad k \in \{1, 2, \dots, \mu^j \}. \label{Eq5}
\end{align}
If we further define the following mapping $K_{\neg j}^0 \colon \mathscr{V}_j \to \mathscr{U}_{\neg j}$
\begin{align}
K_{\neg j}^0 v_k^j = u_k^{\neg j}, \qquad k \in \{1, 2, \dots, \mu^j \}, \label{Eq6} 
\end{align}
 and then by letting $K_{\neg j}$ to be any extension of $K_{\neg j}^0$ to $\mathscr{X}$. We, therefore, have $(A + \sum_{i \in \mathcal{N}_{\neg j}} B_i K_i) v_k^j = w_k^j \in \mathscr{V}$, i.e., $(A + \sum_{i \in \mathcal{N}_{\neg j}} B_i K_i) \mathscr{V}_j \subset \mathscr{V}_j$, so that the controlled-invariant subspace $\mathscr{V}_j$ satisfies 
\begin{align}
 \mathscr{V}_j \in \mathfrak{I}_j (A, B_{\neg j}; \mathscr{X}), \quad \forall j \in \mathcal{N}. \label{Eq7}
\end{align}
\end{IEEEproof}
\begin{corollary} \label{C1}
Let the subspace $\mathscr{V} \subset \mathscr{X}$ be $(A, \, B_{\neg j})$\,-\,invariant for each $j \in \mathcal{N}\cup\{0\}$, then there exists a class of maps $\mathcal{K} \ni K:\mathscr{X} \to \mathscr{U}$ that satisfies
\begin{align}
 \mathcal{K} \subseteq \biggm \{\underbrace{\bigm (K_1, K_2, \ldots, K_N \bigm)}_{\triangleq K} \in \prod_{j \in\mathcal{N}} \mathbb{R}^{r_j \times n} \, \Bigl\lvert \, \biggm(A + \sum_{i \in \mathcal{N}_{\neg j}} B_i K_i\biggm)\mathscr{V} \subset \mathscr{V}, ~~ \forall j \in \mathcal{N} \cup \{0\} \biggm \}. \label{Eq8}
\end{align}
\end{corollary}
\begin{remark}
Note that the controlled-invariant subspace $\mathscr{V}$, which is given in the aforementioned corollary, is also a subspace of $\mathscr{V}^*$ (see also Equation~\eqref{Eq9} below)    
\end{remark}
Next, we introduce the following theorem on the family of supremal controlled-invariant subspaces that will be useful for our work in the next section.
\begin{theorem}\label{T2}
Let $\mathfrak{V}\triangleq\{\mathscr{V}_j^*\}_{j \in \mathcal{N}}$ be a set of supremal controlled-invariant subspaces with respect to the family of systems $\left \{\bigl(A, B_{\neg j} \bigr)\right \}_{j \in \mathcal{N}}$. Then, the set $\mathfrak{V}$ forms a lattice of controlled-invariant subspaces. Moreover, there exists a unique (nonempty) subspace that satisfies 
\begin{align}
 \mathscr{V}^* = \bigcap_{j \in \mathcal{N}}  \mathscr{V}_{j}^* \in \left \{ \bigcap_{j \in \mathcal{N}} \mathscr{V}_j ~ \biggm\vert ~ \mathscr{V}_j \subset \sup \mathfrak{I}_j (A, B_{\neg j}; \mathscr{X}), ~ \forall  j \in \mathcal{N}, ~~ \exists\,u_{\neg j} \in \mathscr{U}_{\neg j} \right \}, \label{Eq9}
\end{align}
where  $\mathscr{V}_{j}^* = \sup \mathfrak{I}_{j} \left ( (A + \sum_{i \in \mathcal{N}_{\neg j}} B_i K_i ), \mathscr{B}_j \right )$ and $\mathscr{B}_j \triangleq \operatorname{Im} B_j$ for all $j \in \mathcal{N}$. 
\end{theorem}
\begin{IEEEproof}
Note that $\mathscr{V}_j + \mathscr{V}_{\neg j} \in \mathfrak{V}$ and $\mathscr{V}_j \cap \mathscr{V}_{\neg j} \in \mathfrak{V}$ for all $j \in \mathcal{N}$. Moreover, if we define the gap metric $\varrho_j(\mathscr{V}_0, \mathscr{V}_j)$ between the controlled-invariant subspaces $\mathscr{V}_0$ and $\mathscr{V}_j$ as
\begin{align}
\varrho_j(\mathscr{V}_0, \mathscr{V}_j) = \Vert P_{\mathscr{V}_0} - P_{\mathscr{V}_j} \Vert, ~~\forall j \in \mathcal{N}, \label{Eq10}
\end{align}
where $\mathscr{V}_0 = \sup \mathfrak{I}_0 (A, B_{\neg 0}; \mathscr{X})$, $P_{\mathscr{V}_0}$ and $P_{\mathscr{V}_j}$ are orthogonal projectors on $\mathscr{V}_0$ and $\mathscr{V}_j$, respectively. Then, the set of all controlled-invariant subspaces in $\mathscr{X}$ forms a compact metric state-space with respect to the above gap metric (see also \cite{GohLR06}). On the other hand, let us define the following family of subspaces
\begin{align}
 \tilde{\mathscr{V}} = \left \{ \bigcap_{j \in \mathcal{N}} \mathscr{V}_j ~ \biggm\vert ~ \mathscr{V}_j \subset \sup \mathfrak{I}_j (A, B_{\neg j}; \mathscr{X}), ~ \forall  j \in \mathcal{N}, ~~ \exists\,u_{\neg j} \in \mathscr{U}_{\neg j} \right \}. \label{Eq11}
\end{align}
Suppose the subspace $\mathscr{V}^*$ exists, then it is a unique member of the family that is defined in \eqref{Eq11}, i.e., 
\begin{align}
 \mathscr{V}^* = \bigcap_{j \in \mathcal{N}}  \mathscr{V}_{j}^* \in \tilde{\mathscr{V}}, \label{Eq12}
\end{align}
with $\mathscr{V}_{j}^* = \sup \mathfrak{I}_{j} \left ( (A + \sum_{i \in \mathcal{N}_{\neg j}} B_i K_i ), \mathscr{B}_j \right )$ for all $j \in \mathcal{N}$. Note that we have $\mathscr{V}_j = \langle A + \sum_{i \in \mathcal{N}_{\neg j}} B_i K_i  \vert \mathscr{B}_j \rangle$ which also implies that $\operatorname{Im} B_j \subset \mathscr{V}_j^*$.\footnote{We remark that the induced continuous maps in $\mathscr{X}/\mathscr{V}_j^*$ and $\mathscr{X}/\mathscr{V}_{\neg j}^*$ admit an enveloping lattice for the family of controlled-invariant subspaces $\mathscr{V}_j^*$, $\forall j \in \mathcal{N}$ (e.g., see \cite{GohLR06}).}
\end{IEEEproof}

\subsection{Properties of (controlled)-invariance entropy}
In the following, we start by giving the definition of (controlled)-invariance entropy for the multi-channel system in \eqref{Eq1} with respect to the subset of control-input channels and that class of control functions. 
\begin{definition}
For a given subspace $\mathscr{F} \subset \tilde{\mathscr{V}}^* \in \tilde{\mathscr{V}}$ with nonempty interior and $T$, $\epsilon > 0$, the class of control functions  $\mathscr{C}(T, \epsilon, \mathscr{F}, \tilde{\mathscr{V}}^*)\subset \mathscr{U}$ is called $(T, \epsilon, \mathscr{F}, \tilde{\mathscr{V}}^*)$-spanning, if there exits $u \in \mathscr{C}(T, \epsilon, \mathscr{F}, \tilde{\mathscr{V}}^*)$ {\em for almost all} $t \in [0,\, T)$ such that
\begin{align}
 \max_{\substack{ j \in \mathcal{N}}} \sup_{\substack{t \in [0,\, T]}} \inf_{\substack{ y \in \tilde{\mathscr{V}}^*}} \Vert \phi_{\neg j} (t, x_0, u_{\neg j}(t)) - y\Vert \le \epsilon, \quad \forall x_0 \in \mathscr{F}. \label{Eq13}
\end{align}
\end{definition}
\begin{remark}
In the aforementioned definition, we use the notation $\phi_{\neg j}(t, x_0, u_{\neg j}(t))$ to denote the unique solution of the multi-channel system with initial condition $x_0 \in \mathscr{F}$ and control $u_{\neg j} \in \mathscr{U}_{\neg j}$, i.e., 
\begin{align}
 x(t) &= \phi_{\neg j}(t, x_0, u_{\neg j}(t)), \notag \\
       &\triangleq \exp A \bigl(t - t_0 \bigr)x_0 + \sum_{i \in \mathcal{N}_{\neg j}} \int_{t_0}^t \exp A\bigl(t - s \bigr) B_i u_i(s)ds, \quad \forall [t_0,\,t]\in[0,\,T],\label{Eq14}
\end{align}
for each $j \in \mathcal{N}\cup\{0\}$. 

Moreover, the relation $\phi_{\neg j}(t + t_0, x_0, u_{\neg j}(t)) = \phi_{\neg j}(t, \phi_{\neg j}(t, x_0, u_{\neg j}(t)), u_{\neg j}(t_0 + .))$ will also hold for all $j \in \mathcal{N}\cup\{0\}$.
\end{remark}

Let $r_{\rm inv}(T, \epsilon, \mathscr{F}, \tilde{\mathscr{V}}^*)$ be the smallest cardinality of $\mathscr{C}(T, \epsilon, \mathscr{F}, \tilde{\mathscr{V}}^*)$-spanning sets. Then, we have the following properties for $r_{\rm inv}(T, \epsilon, \mathscr{F}, \tilde{\mathscr{V}}^*)$.
\begin{enumerate} [(i)]
\item Clearly $r_{\rm inv}(T, \epsilon, \mathscr{F}, \tilde{\mathscr{V}}^*) \in [0,\,\infty)$.\footnote{The value of $r_{\rm inv}(T, \epsilon, \mathscr{F}, \tilde{\mathscr{V}}^*)$ could be an infinity.}
\item If $\epsilon_1 < \epsilon_2$, then $r_{\rm inv}(T, \epsilon_1, \mathscr{F}, \tilde{\mathscr{V}}^*) \ge r_{\rm inv}(T, \epsilon_2, \mathscr{F}, \tilde{\mathscr{V}}^*)$.
\end{enumerate}

\begin{definition}
The (controlled)-invariance entropy of the multi-channel system in \eqref{Eq1} (i.e., with respect to the subset of control-input channels and/or the class of control functions) is given by 
\begin{align}
 {h_{\rm inv}}(\mathscr{F}, \tilde{\mathscr{V}}^*) = \lim_{\substack{\epsilon \searrow 0}} \biggm\{\varlimsup_{\substack{T \to \infty}} \frac{1}{T} \log r_{\rm inv}(T, \epsilon, \mathscr{F}, \tilde{\mathscr{V}}^*)\biggm\}. \label{Eq15}
\end{align}
\end{definition}

\begin{remark}
We remark that the existence of such a limit in the aforementioned definition for ${h_{\rm inv}}(\mathscr{F}, \tilde{\mathscr{V}}^*)$ follows directly from the monotonicity of $r_{\rm inv}(T, \epsilon, \mathscr{F}, \tilde{\mathscr{V}}^*)$ with respect to $\epsilon$. Moreover, such an invariance entropy ${h_{\rm inv}}(\mathscr{F}, \tilde{\mathscr{V}}^*)$ equals to the minimum amount of information that is required to render $\tilde{\mathscr{V}}^*$-invariant subspace by using a causal coding and/or control law (see \cite{Colo10} for discussion on single control-channel systems).
\end{remark}

Then, we have the following properties for ${h_{\rm inv}}(\mathscr{F}, \tilde{\mathscr{V}}^*)$.
\begin{enumerate} [(i)]
\item ${h_{\rm inv}}(\mathscr{F}, \tilde{\mathscr{V}}^*) \in [0,\,\infty) \cup\{\infty\}$.
\item If $\mathscr{F}\triangleq \bigcup_{l \in \{1,2, \ldots, L\}} \mathscr{F}_l$ with compact $\mathscr{F}_l$, $\forall l \in \{1,2, \ldots, L\}$, then ~ ${h_{\rm inv}}(\mathscr{F}, \tilde{\mathscr{V}}^*)  =  \max_{\substack{l \in \{1,2, \ldots, L\}}} {h_{\rm inv}}(\mathscr{F}_l, \tilde{\mathscr{V}}^*)$.
\end{enumerate}

In the following, we state the main problem of this paper -- where we establish a connection between the invariance entropy of a multi-channel system and the reliability of a multi-channel system.

{\em Problem}: Find a condition on the minimum amount of ``information" (with respect to the subset of control-input channels and/or the class of control functions) that is necessary to keep the states of the multi-channel system in a given subspace $\tilde{\mathscr{V}}^*$.

\section{Main Results} \label{sec:main}
In this section, we present our main results -- where we provide a connection between the invariance entropy (as a measure of ``information" needed with respect to the subset of control-input channels to keep the system in (or near) this compact subspace) and the reliability of the multi-channel system (when there is an intermittent channel failure that may occur between actuators, controllers and/or sensors in the system).
\begin{theorem} \label{T3}
Suppose that Theorem~\ref{T2} holds true and let $\mathscr{F}$ be a subspace of $\tilde{\mathscr{V}}^*$.\, For every $x_0 \in \mathscr{F}$, if there exists $u(t) \in \mathscr{C}(T, \epsilon, \mathscr{F}, \tilde{\mathscr{V}}^*)$ {\em for almost all} $t \in [0,\, T]$ that renders $\tilde{\mathscr{V}}^*$-invariant, then the (controlled)-invariance entropy of the multi-channel system is given by 
\begin{align}
 h_{\rm inv}(\mathscr{F}, \tilde{\mathscr{V}}^*) = \lim_{\substack{\epsilon \searrow 0}}\biggm\{ \varlimsup_{\substack{T \to \infty}} \frac{1}{T} \log r_{\rm inv}(T, \epsilon, \mathscr{F}, \tilde{\mathscr{V}}^*)\biggm\}. \label{Eq16}
\end{align}
\end{theorem}
\begin{IEEEproof}
For any subspace $\mathscr{F} \subset \tilde{\mathscr{V}}^*$, suppose there exists $u(t) \in \mathscr{C}(T, \epsilon, \mathscr{F}, \tilde{\mathscr{V}}^*)$ {\em for almost all} $t \in [0,\, T]$ such that  
\begin{align}
\sup_{\substack{ t \in [0,\,T]}} \inf_{\substack{ y \in \tilde{\mathscr{V}}^*}} \Vert \phi_{\neg j} (t, x_0, u_{\neg j}(t)) - y\Vert \le \epsilon, ~~\forall j \in \mathcal{N}, ~~  \forall x_0 \in \mathscr{F}, \notag
\end{align}
with some finite-positive number $\epsilon>0$. Then, we see that \eqref{Eq16} (i.e., the (controlled)-invariance entropy of the multi-channel system in \eqref{Eq1}) will directly follow. Moreover, the supremum is taken overall the set of all admissible controls that renders $\tilde{\mathscr{V}}^*$-invariant.
\end{IEEEproof}
\begin{remark}
We remark that the aforementioned theorem essentially states that the set of admissible controls $\mathscr{U}$ renders $\tilde{\mathscr{V}}^*$-invariant, even if there is an intermittent failure in any one of the control-input channels. Note that this quantity (which is also the minimum growth rate for the number of open-loop control functions with respect to intermittently faulty channel) provides a condition on the minimum amount of ``information" that is necessary to keep the system states in (or near) this subspace.
\end{remark}
We conclude this section with the following result which is an immediate corollary of Theorem~\ref{T3}.
\begin{corollary} \label{C2} 
Suppose there exists a finite-positive number $\epsilon_{\mathcal{K}} > 0$ such that
\begin{align}
\max_{\substack{ j \in \mathcal{N}}} \sup_{\substack{ t \in [0,\,T]}} \inf_{\substack{ y \in \tilde{\mathscr{V}}_{\mathcal{K}}^*}} \Vert \phi_{\neg j} (t, x_0, u_{\neg j}(t)) - y\Vert \le \epsilon, \quad \forall x_0 \in \mathscr{F}, \label{Eq17}
\end{align}
where
\begin{align}
 \tilde{\mathscr{V}}_{\mathcal{K}}^* = \sup \left \{ \bigcap_{j \in \mathcal{N}} \mathscr{V}_j ~ \biggm\vert ~ \mathscr{V}_j \subset \sup \mathfrak{I}_j (A, B_{\neg j}; \mathscr{X}), ~ \forall  j \in \mathcal{N}, ~~ \exists\,K \in \mathcal{K} \right \} \supset \mathscr{F}. \label{Eq18}
\end{align}
Then, the (controlled)-invariance entropy of the multi-channel system in \eqref{Eq1} is given by 
\begin{align}
 h_{\rm inv}(\mathscr{F}, \tilde{\mathscr{V}}_{\mathcal{K}}^*) = \lim_{\substack{\epsilon_{\mathcal{K}} \searrow 0}} \biggm\{\varlimsup_{\substack{T \to \infty}} \frac{1}{T} \log r_{\rm inv}(T, \epsilon_{\mathcal{K}}, \mathscr{F}, \tilde{\mathscr{V}}_{\mathcal{K}}^*)\biggm\}. \label{Eq19}
\end{align}
\end{corollary}
\begin{remark}
Note that the bounds for $h_{\rm inv}(\mathscr{F}, \tilde{\mathscr{V}}^*)$ and $h_{\rm inv}(\mathscr{F}, \tilde{\mathscr{V}}_{\mathcal{K}}^*)$ are different, since they may depend on their respective classes of control functions. Moreover, the following $h_{\rm inv}(\mathscr{F}, \tilde{\mathscr{V}}^*) \le h_{\rm inv}(\mathscr{F}, \tilde{\mathscr{V}}_{\mathcal{K}}^*)$ also holds true.
\end{remark}

\end{document}